\documentclass[11pt]{amsart}
\usepackage{latexsym,amssymb,amsmath}
\def\boxit#1{\vbox{\hrule height1pt\hbox{\vrule width1pt\kern3pt
  \vbox{\kern3pt#1\kern3pt}\kern3pt\vrule width1pt}\hrule height1pt}}

\def\qq#1#2#3{q^{#1}_{{#2} {#3}}}
\def\ii{|II_{X,x}|}
\def\rr#1#2#3#4{r^{#1}_{{#2} {#3}{#4}}}
\def\rrr#1#2#3#4#5{r^{#1}_{{#2} {#3}{#4}{#5}}}

\def\tmod{\text{ mod }}
\def\tbase{\text{Base}}
\def\tsing{\text{Sing}}
\def\ooo#1#2{\omega^{#1}_{#2}}
\def\oo#1{\omega^{#1}_0}
\def\ee#1{e_{#1}}
\def\ue#1{e^{#1}}
\def\ep{\epsilon}
\def\cf{\mathcal F}
\def\xx#1{x^{#1}}
\def\frp#1#2{\frac{\partial {#1}}{\partial {#2}}}

\def\BC{\mathbb C}\def\BO{\mathbb O}\def\BS{\mathbb S}
\def\BA{\mathbb A}\def\BR{\mathbb R}\def\BH{\mathbb H}
\def\BP{\mathbb P}
\def\pp#1{\mathbb P^{#1}}
\def\fr{\mathfrak r}

\def\fd{\mathfrak d}

\def\pp#1{{\mathbb P}^{#1}}
\def\tdim{\rm dim}
\def\hd{,...,}
\def\ww{\wedge}

\def\na{n+a}

\def\11{\mathbf 1}

\def\fsl{{\mathfrak {sl}}}

\def\fg{{\mathfrak g}}

\def\a{\alpha}

\def\o{\omega}

\def\b{\beta}
\def\g{\gamma}

\def\ot{{\mathord{\,\otimes }\,}}
\def\op{{\mathord{\,\oplus }\,}}

\def\ra{{\mathord{\;\rightarrow\;}}}

\def\La#1{\Lambda^{#1}}

\newtheorem{theorem}{Theorem}[section]

\theoremstyle{definition}

\newtheorem{example}[theorem]{Example}

\theoremstyle{remark}
\newtheorem{remark}[theorem]{Remark}

\begin{document}

\title{Griffiths-Harris rigidity of compact Hermitian symmetric spaces}
\author{J.M. Landsberg}
\begin{abstract} I show that any complex manifold that resembles
a rank two compact Hermitian symmetric space (other than a quadric
hypersurface) to order two at a general point must be an open
subset of such a space.
\end{abstract}
\maketitle


\section{Introduction}

Let $X\subset\BC\pp\na$ be a variety and let $x\in X $ be a smooth point.
The projective second fundamental form of $X$ at $x$ (see 
\cite{GH,Lsnu,Lrr,EDSBK})
is a basic differential invariant that measures how $X$ is moving
away from its embedded tangent projective space at $x$ to first order.
It determines a system of quadrics $\ii\subset S^2T^*_xX$.
Now let $X$ be a general point, {\it to what extent does
$\ii$ determine $X$?}

Let $X$ be such that $\ii$ is an isolated
 point in the moduli space of $a$-dimensional
linear subspaces of the space  of quadratic forms on $\BC^n$
up to linear equivalance. We say $X$ is {\it infinitesimally rigid at 
order two} or is {\it Griffiths-Harris rigid} if whenever 
  $Y\subset \pp N$ is a complex manifold,
  $y\in Y$ is a general point  and
$|II_{Y,y}|=\ii$, then $\overline Y=X$.

In \cite{GH}, Griffiths and Harris posed the question as to whether the 
Segre variety $Seg(\pp 2\times \pp 2)\subset \pp 8$ was infinitesimally
rigid to order two and in \cite{Lrigid} I 
answered the question affirmatively and showed that  all rank two compact Hermitian
symmetric spaces (in their minimal homogeneous embeddings)
except for the quadric hypersurface, the Grassmanian $G(\BC^2,\BC^5)\subset \pp 9$ and the {\it spinor
variety}  $D_5/P_5=\BS_{10}\subset\pp{15}$ were
infinitesimally rigid at order two.
The quadric is not rigid to order two and
Fubini showed \cite{fub} it is rigid to order three when
$n>1$ and   it is rigid to order five when $n=1$.   In this paper I
  resolve the two remaining cases,
and explain   shorter  and less computational proofs for the other
cases presented in \cite{Lrigid}. I also reprove the
rigidity of the three Severi varieties that are rigid
to order two to illustrate the method.  The new proofs use two tools,
a higher order Bertini theorem, and elementary representation theory.

In \cite{LM0, LM1} we showed that  all   rational homogeneous
varieties other than the rank two compact Hermitian symmetric spaces
fail to be rigid to order two, so the result of this paper is the
best possible in this sense. One can compare this type of rigidity
to that studied by Hwang and Mok, see, \cite{Hwang, HM}.
Some differences are: in their study they require  global
hypotheses where here the hypotheses are at the level of
germs (this is because the systems of quadrics under study admit
no local deformations); in their study the objects of interest are
not  {\it a priori} given an embedding (although since they assume the Picard
number is one, one gets something close to an embedding); and in their
study the object of interest is the cone of minimal degree rational
curves through a general point, which, {\it a priori}, has nothing
to do with the cone of asymptotic directions I use here (in the systems under
consideration, the base locus of $II$ determines $II$).

\begin{theorem} Let $X^n\subset \BC\pp\na$ be a complex submanifold.
Let $x\in X$ be a general point. If $\ii \simeq |II_{Z,z}|$ where $Z$ is
a   compact rank two Hermitian symmetric space in its
natural embedding, other than
a quadric hypersurface,  then $\overline X =Z$.\end{theorem}

Some open questions and relations with the Fulton-Hansen connectedness
theorem are discussed in \cite{Lrigid}. Another application of the
  techniques used here is given
in \cite{LMleg}.

\section{Moving frames}

For more details throughout this section, see any of   \cite{GH,Lci,Lsnu,EDSBK}.

Once and for all fix
 index ranges $1\leq \a,\b,\g\leq n$, $n+1\leq \mu,\nu\leq\na$.

Let $X^n\subset\BC\pp\na=\BP V$ be a 
complex submanifold and let $x\in X$ be a general
point. Let $\pi :\cf^1\ra X$ denote the bundle of bases of $V$ (frames) preserving
the flag
$$
\hat x\subset\hat T_xX\subset V.
$$
Here $\hat T_xX$ denotes the affine tangent space (the cone over the
embedded tangent projective space).
Let $(\ee 0\hd \ee\na)$ be a basis of $V$ with dual basis
$(\ue 0\hd \ue\na)$
  adapted such that $\ee 0\in \hat x$ and
$\{\ee 0,\ee \a\}$ span $\hat T_xX$.
I ignore twists and obvious quotients, 
writing $\ee\a$ for $(\ee \a \tmod \ee 0)\ot e^0\in T_xX$ and
$\ee\mu$ for $(\ee \mu \tmod \hat T_xX )\ot e^0\in N_xX=T_x\BP V/T_xX$.
Moreover, if $x$ and $X$ understood, I write $T=T_xX$ and
$N=N_xX$.

The fiber of $\pi : \cf^1\ra X$ over
a point is isomorphic to the group
$$
G_1=\left\{
\ g = \begin{pmatrix} g^0_0 & g^0_{\beta} & g^0_{\nu} \\
0 &  g^{\alpha}_{\beta} & g^{\alpha}_{\nu} \\
0 &  0 & g^{\mu}_{\nu} \end{pmatrix} \mid g\in GL(V)  \right\}.
$$

While $\cf^1$ is not in general a Lie group, since $\cf^1\subset GL(V)$,
we may pullback the Maurer-Cartan from on $GL(V)$ to $\cf^1$.
Write the pullback  of the Maurer-Cartan form to
$\cf^1$ as
$$
\omega=\begin{pmatrix}\oo 0 & \ooo 0\beta & \ooo 0\nu\\
\oo\alpha & \ooo\alpha\beta & \ooo\alpha\nu\\
\ooo\mu 0 & \ooo\mu\beta & \ooo\mu\nu\end{pmatrix} .
$$
The adaptation implies that $\oo\mu =0$ and
the Maurer-Cartan equation $d\omega =-\omega\ww\omega$ together
with the Cartan Lemma implies that
 for all
$\mu,\a$,  $\ooo\mu\a=\qq\mu\a\b\oo\b$ for some
functions $\qq\mu\a\b=\qq\mu\b\a :\cf^1\ra\BC$. These functions determine the
projective second fundamental form $II=F_2=\ooo\mu\a\ot\ee\mu=
\qq\mu\a\b\oo\a\oo\b\ot \ee\mu \in \Gamma (X,S^2T^*X\ot NX)$.

While $II$ descends to be a section of $(S^2T^*X\ot NX)$, 
higher order derivatives provide   relative differential invariants
 $F_k\in \Gamma (\cf^1, \pi^*(S^kT^*\ot N))$. For example,
$$
\begin{aligned}
F_3&= \rr\mu\alpha\beta\gamma\oo\a\oo\b\oo\gamma\ot\ee\mu \label{f3coeff}\\
F_4 &=  \rr\mu\alpha\beta{\gamma\delta}\oo\delta  \oo\a \oo\b \oo\g \ot \ee\mu 
\end{aligned}
$$
where the functions $\rr\mu\alpha\beta\gamma,\rr\mu\alpha\beta{\gamma\delta}$
are given by

\begin{equation}
\rr\mu\alpha\beta\gamma\oo\gamma =
-d\qq\mu\alpha\beta - \qq\mu\alpha\beta\ooo 00 -\qq\nu\alpha\beta\ooo\mu\nu
+\qq\mu\alpha\delta\ooo\delta\beta +
 \qq\mu\beta\delta\ooo\delta\alpha \label{f3eqn}
\end{equation}
\begin{equation}
\rr\mu\alpha\beta{\gamma\delta}\oo\delta  
 =  
-d\rr\mu\alpha\beta\gamma -2\rr\mu\alpha\beta\gamma\ooo 00
-\rr\nu\alpha\beta\gamma\ooo\mu\nu +
\mathfrak  S_{\alpha\beta\gamma}(
\rr\mu\alpha\beta\epsilon\ooo\epsilon\gamma +3\qq\mu\alpha\beta\ooo 0\gamma
-\qq\mu\alpha\epsilon\qq\nu\beta\gamma\ooo\epsilon\nu ). \label{f4eqn}
\end{equation}

If one chooses local affine coordinates and writes $X$ as a graph
$$
\xx\mu = \qq\mu\a\b\xx\a\xx\b + \rr\mu\a\b\g\xx\a\xx\b\xx\g
+\rrr\mu\a\b\g\delta \xx\a\xx\b\xx\g\xx\delta + ...
$$
then there exists a  local section of $\cf^1$
such that 
$$
\begin{aligned}
F_2&=\qq\mu\a\b d\xx\a d\xx\b\ot\frp{}{\xx\mu}\\
F_3&=\rr\mu\a\b\g d\xx\a d\xx\b d\xx\g\ot\frp{}{\xx\mu}\\
F_4&=\rrr\mu\a\b\g\delta d\xx\a d\xx\b d\xx\g d\xx\delta\ot\frp{}{\xx\mu}
\end{aligned}
$$
etc...

Since an analytic variety is uniquely determined by its Taylor series
at a point, to show $Z$ is rigid to order
two,  it is sufficient to show that over varieties $X$ with
$|II_{X,x}|=|II_{Z,z}|$ there exists a subbundle of $\cf^1$ such
that the $F_k$'s of $X$ coincide with those of $Z$. Moreover,
the minuscule varieties, that is, the compact Hermitian symmetric
spaces in their natural projective embeddings,
 have the property that  on a reduced frame
bundle all the differential invariants except for their
fundamental forms are zero and in our case the only nonzero fundamental
form is $II$.
 
The method in \cite{Lrigid} was first to use the equations above to calculate
relations among the coefficients of $F_3$. Enough relations were
found that, combined with the coefficients that were normalizable to zero,
I obtained that $F_3$ was zero, and the same technique was
used for higher order invariants.

In this paper I   decompose $S^3T^*\ot N$ into irreducible $R$-modules,
where $R\subset GL(T)\times GL(N)$ is the subgroup preserving
$II\in S^2T^*\ot N$. I also systematize the vanishing of coefficients of the $F_k$ somewhat
 using higher order Bertini theorems which I now describe. It would be
nice to have a way to apply the higher order Bertini theorems directly
to the irreducible modules instead of using the coefficients of $F_3$.

\section{Vanishing tools}

\subsection{Higher order Bertini}
Let $T$ be a vector space.
The classical Bertini theorem implies that for
a linear subspace $A\subset S^2T^*$,
  if $q\in A$ is generic,
then $v\in Sing(q):=\{ v\in T\mid q(v,w)=0\ \forall w\in T\}$ 
implies $v\in Base(A):=\{ v\in T\mid Q(v,v)=0\ \forall Q\in A\}$.

\begin{theorem}[Higher order Bertini] Let $X^n\subset \BP V$ be
a complex manifold and let $x\in X$ be a general point.  

1.  Let $q\in\ii$ be a generic quadric. 
Then $q_{sing}\subset\tbase \{ F_2\hd F_k\}$ for all $k$. I.e.,
$q_{sing}$ is tangent to a linear space on  the completion of  $X$.
  
2. Let $q\in\ii$ be a any quadric, let $L\subset q_{sing}\cap\tbase\ii$
  be a linear subspace.
Then for all $v,w\in L$, $F^q_3(v,w,\cdot)=0$. 

3. With $L$ as in 2.,
if $L'\subset (\tbase\{\ii, F_3\})\cap L$ is a linear subspace then 
$F^q_4(u,v,w,\cdot)=0$ for all $u,v,w\in L'$ and so on for higher orders.
Here $F^q_k$ denotes the polynomial in $F_k$ corresponding to the
conormal direction of $q$. This is well defined by the lower order vanishing.

4. With $L'$ as in 3., if $L''\subset L'\cap (F_3^q)_{sing}$ is a linear
space,
then for all $u,v\in L''$, $F^q_4(u,v,\cdot,\cdot)=0$.

\end{theorem}

Analogous results hold for higher orders.

\begin{proof} Note that 1. is classical, but we provide a proof
for completeness.
Assume $v=\underline e_1$ and $q=q^{\mu}$. Our hypotheses
imply $\qq\mu 1\b=0$ for all $\b$. Formula \eqref{f3eqn} reduces to  
$$
\rr\mu 11\b\oo\b= -\qq\nu 11\ooo\mu\nu .
$$
If $q$ is generic we are still working on $\cf^1$
and so the $\ooo\mu\nu$ are independent of
the semi-basic forms, thus the coefficients on both sides of the equality
are zero, proving both the classical
Bertini theorem and 1 in the case $k=3$.  
 
 If $q$ is not generic, in order to have $q=q^{\mu},v=\ee 1$ we have reduced
 to a subbundle $\cf'\subset \cf^1$ and we
no longer have the  $\ooo\mu\nu$ independent.
However hypothesis 2 states that $\qq\nu 11=0$ for all $\nu$ and the 
required vanishing still holds. 

For 3., note that
$\rrr\mu 111\delta\oo\delta =\rr\nu 111\ooo\mu\nu + e\rr\mu 11\ep \oo\ep 1
+\qq\mu 1\ep \qq\nu 11$ and the right hand side is zero under
our hypotheses. Part 4 is proven simliarly.  

The extension to linear spaces holds by polarizing the forms.
The analogous equation at each order proves the next higher order.
\end{proof}

 \begin{example} Let $X=G(2,m)$ and let $V=\La 2\BC^m$ have basis $e_{st}$ with $1\leq s<t\leq m$.
At $x=[e_{12}]$ we have the adapted flag
 $$
\{e_{12}\} \subset \{e_{1j},e_{2j}\} \subset V
$$
where $3\leq i<j\leq m$,
and $SL_2\times SL_{m-2}$ acts transitively on $N_x\simeq\{ e_{ij} \}$.
So here $\a=\{ (1j),(2,j)\}$, $\mu = \{ (ij) \}$.
In these frames $II=(\oo{(1i)}\oo{(2j)}-\oo{(1j)}\oo{(2i)})\ot\ee{ij}$.

If $m=5$, then 
  $q^{45}$ is a generic quadric with $e_{13}\subset q^{45}_{sing}$.
Thus we have
$$
\begin{aligned}
\rr{\mu}{(13)}{(13)}{(13)} &=0 \ \forall \mu\\
\rr{45}{(13)}{(13)}{\b} &=0\ \forall\b .\end{aligned}
$$

If $m>5$ then $q^{45}$ is no longer generic, but since
$\ee{(13)}\in\tbase\ii$ we still may conclude
$$
\rr{45}{(13)}{(13)}{\b}  =0\  \forall \b .
$$
\end{example}

\subsection{Normalizations} $F_3$ is translated in the fiber of $\cf^1$ by the
action of $T\ot N^*$ and $T^*$ (the $g^{\a}_{\mu}$ and the $g^0_{\a}$).
We may decompose $T\ot N^*$ and $T^*$ into irreducible $R$ modules
and determine which of these act nontrivially. In the case the variety is modeled
on a
rank two minuscule variety, we will have that all of $T\ot N^*$ acts
effectively, but the $T^*$ action duplicates a factor in $T\ot N^*$.
This is because in the homogeneous model, the forms $\ooo 0\b$
are independent and the forms $\ooo\a\mu$ are linear combinations
of the $\ooo 0\b$. We will let $\cf^n$ denote the bundle where
the action of $T\ot N^*$ has been used to kill the corresponding components
of $F_3$. Similarly, on $\cf^n$, $F_4$ is translated by the action
of $N$ and we will let $\cf^N$ denote the subbundle of $\cf^n$ where
the component of $N$ in $F_4$ has been normalized to zero.

\subsection{Remarks on   decompositions of the $F_k$ and vanishing}
  Let $II\in S^2T^*\ot N$ arise from a trivial representation
of a  reductive  group $R\subset GL(T)\times GL(N)$.
Let $X^n\subset\pp\na$ be a complex submanifold, let $x\in X$ be a general
point and suppose
$II_{X,x}=II$.  
\smallskip

1. Since the oribit of a highest weight vector in any module spans
the module, the component of $F_k$ in an irreducible module $V$
is zero if its   highest weight vector is zero.

\smallskip

2. An irreducible module in $S^3T^*\ot N$ can occur in $F_3$ only
if it also occurs in
$(T\ot T^*)^{\fr^c}\ot T^* + (N\ot N^*)^{\fr^c}\ot T^*$.
Here $\fr$, the Lie algebra of $R$, occurs as a submodule of
$T\ot T^*$ and $N\ot N^*$ and  $(T\ot T^*)^{\fr^c}$ denotes the
complement of $\fr$ in $T\ot T^*$.
This is because  the tangential and normal
connection forms, $\ooo\a\b,\ooo\mu\nu$
may be decomposed into $\rho_T(\fr)$-valued (resp. $\rho_N(\fr)$-valued)
forms and
semi-basic forms.
(Here $\rho_T,\rho_N$ denote the representations of $\fr$ on
$T$ and $N$.)   The coefficients for the semi-basic components are linear combinations
of  the $\rr\mu\a\b\g$, 
all of the same weight. On the other hand the coefficients give
rise to $R$-modules respectively in 
$(T\ot T^*)^{\fr^c}\ot T^*$ and $(N\ot N^*)^{\fr^c}\ot T^*$.

\smallskip

3. Similarly, if $F_3=0$, and the normalizations of $F_3$
are exactly by $T \ot N^*$, then an irreducible 
module in $S^4T^*\ot N$ can occur in $F_4$ only
if it also occurs in
$(T\ot N^*)^{T^{*c}}\ot T^*$,  where $(T\ot N^*)^{T^{*c}} $
is the complement of $T^*$ in $T\ot N^*$. Thus in our normalizations, the
forms $\ooo 0\b$ will remain independent and independent of the semi-basic
forms. On the other hand the forms $\ooo\ep\nu$ will become dependent
on the semi-basic forms and the $\ooo 0\b$. Again, the components that
will depend on the semi-basic forms will have coefficients consisting
of linear combinations of monomials in $F_4$ of the same weight.

\smallskip

4. If $F_3,F_4=0$ (after normalizations), then an irreducible 
module in $S^5T^*\ot N$ can occur in $F_4$ only
if it also occurs in $N$.

\smallskip

5. If $F_3,F_4,F_5=0$ (after normalizations), then all higher
$F_k$ are zero as well, see \cite{Lci}.

\smallskip

\section{Case of $G(2,5)$ and $\BS_{10}$}

\subsection{Model for $G(2,5)$}
Write $T=  A^*\ot B$. We index   bases of $T$ and $N$ as above.
  $R =   \fsl(A)+\fsl(B)+\BC=\fsl_2+\fsl_3+\BC$. We write $A_j$ for the represention
of $\fsl(A)$ with highest weight $j$ and $B_{ij}$ for the representation
of $\fsl(B)$ of highest weight $i\omega_1+j\omega_2$. Here and throughout
we use the notations and ordering of the weights of \cite{bou}.
The relevant modules are summarized in the table below.

\subsection{Model for $\BS_{10}$}
 Write
 $\BC^{16}=Clifford(\BC^5)\simeq \Lambda^{even}\BC^5$
 with $\hat x\simeq \La 0W$, $T\simeq \La 2W, N\simeq \La 4W\simeq W^*$.
We let $e_{st}=e_s\ww e_t$, $1\leq s<t\leq 5$ index a basis of $T$ and
$e^s$ index a basis of $N$. Note that, as with $G(2,5)$, $R$ acts
transitively on $N$ so all quadrics in $|II|$ are generic.

Let $V_{ijkl}$ denote the $\fsl_5$ module with highest weight $i\o_1+
j\o_2+k\o_3+l\o_4$. $|II|$ is
given by the Pfaffians of the $4\times 4$ minors centered about
the diagonal with  $e^j$ corresponding to the Pfaffian obtained
by removing the $j$-th row and column.
The relevant modules are summarized in the following
table.

$$
\begin{aligned}
T&=V_{0100}\\
&=A_1\ot B_{10}\\
T^*&=V_{0010}\\
&=A_1\ot B_{01}\\
N&= V_{0001 }\\
&=A_0\ot B_{10}\\
S^2T^*&= T^{*2}\op N^*\\
S^3T^*&= T^{*3}\op N^*T^*\\
S^3T^*\ot N&=(T^{*3}N\op TT^{*2})\op ((N^*T^*)N\op N^*T\op T^*)\\
T \ot N^*&=  TN^*\op T^*\\
(T \ot N^*)^{T^{*c}}\ot T^*&=   N^*TT^*\op N^{*2}N\op TN\op N^*\\
 S^4T^*&= T^{*4}\op N^*T^{*2}\op N^{*2}\\
  S^4T^*\ot N&= (T^{*5}\op TT^{*3}N\op N^*T^{*3}) \op (N^*T^{*3}  \op N^*TT^*N \op T^{*2}N \op N^{*2}T^*
  \op TT^*)\\
  &
  \op (N^{*2}N\op N^*)
\end{aligned}
$$
The notation is such that if $V_{\lambda}, W_{\mu}$ are the
irreducible representations with highest weights $\lambda,\mu$, then
$V^{k}$, $VW$ are respectively the representations with highest
weights $k\lambda$ and $\lambda +\mu$. $VW\subset V\ot W$ is
called the {\it Cartan component}.

To obtain the vanishing of $F_3$ we need to eliminate five modules.
We first eliminate two by reducing to $\cf^n$
as described above, so the
last two factors are zero. 
   Let $\cf^n\subset \cf^1$ denote
our new frame bundle.  

On our new bundle there remains three modules to eliminate.

The first module in $S^3T^*$
is decomposably generated by $r_{(13)(13)(13)}$ in the $G(2,5)$
case and $r_{(12)(12)(12)}$ in the $\BS_{10}$ case.
We already saw the $G(2,5)$ case is covered by Bertini,
and the $\BS_{10}$ case is as well because
$\ee{(12)}\in q^{1}_{sing}$, and  all quadrics in the system
are generic. Thus the first two modules in $S^3T^*\ot N$   don't appear in $F_3$.

At this point just $(N^*T^*)N$ remains.   In the $G(2,5)$ case
$(N^*T^*)$ has highest weight a linear combination of
$r_{(13)(13)(24)}$ and $r_{(13)(14)(23)}$. 
In the $\BS_{10}$ case it has highest weight a linear combination
of
$r_{(12)(12)(34)}$,$r_{(12)(13)(24)}$  and $r_{(12)(14)(23)}$ and thus
the Cartan components respectively  have highest weights linear combinations of
$r_{(13)(13)(24)}^{45}$, $r_{(13)(14)(23)}^{45}$
and
$r_{(12)(12)(34)}^5$, $r_{(12)(13)(24)}^5$ and $r_{(12)(14)(23)}^5$, all
of which are zero by Bertini.
Thus $F_3$ is zero.

We   normalize away the $N$ factor in $S^4T^*\ot N$ and
study the remaining modules. Comparing $S^4T^*\ot N$ and
$(T \ot N^*)^{T^{*c}}\ot T^*$ modulo $N^*$, their intersection is
empty and thus $F_4=0$ on $\cf^N$.

One can check that $S^5T^*\ot N$ does not contain a copy of $N$, so
  we are done.\qed

\begin{remark} If one compares modules, $N^*T^*N$ is not eliminated
from $F_3$. On the other hand, if one just uses Bertini to study $F_4$,
  everthing coming from the first factor  in $S^4T^*$
and all Cartan
products of $N$ with factors in $S^4T^*$
are easily seen to be zero, but a few of the other
modules are more complicated to eliminate. \end{remark}

\section{Case of $\BA\pp 2$}

Let $\BA_{\BR}$ respectively denote $\BC,\BH,\BO$ and let
$\BA=\BA_{\BR}\ot_{\BR}\BC$.

I use the following model: $T=\BA\op\BA$, where $\BA=\BH$,
the complexified quaternions, for $G(2,6)$, and $\BA=\BO$, the
complexified octonions for $\BO\pp 2$.
I use
$(a,b)$ as $\BA$-valued coordinates. Then $|II|=\{ a\overline a,
b\overline b, a\overline b\}$ where $a\overline b$ represents
$\tdim\BA$ quadrics.
Let $p=3,7$. Write $a=a_0+a_1\ep_1+\hdots + a_p\ep_p$.
We will need to work with null vectors
 so let $\ee 1=1+i\ep_1$, $\overline e_1=1-i\ep_1$,
$\ee 2=1+i\ep_2$  denote elements of the first copy of $\BA$ (with
coordinate $a$).
We let $\ee a$ denote the normal vector such that $q^a=a\overline a$
and similarly for $\ee b$. Let $\ee 0$ denote the normal vector
such that $q^0=Re(a\overline b)$ and  $\ee {\ep_j}$
such that $q^{\ep_j}$ is the $\ep_j$ coefficient of $a\overline b$.

  Let
$V_{ijklm}$ denote the $\fd_5$-module with highest weight
$i\o_1+j\o_j+k\o_3+l\o_4+m\o_5$, and the $\fsl(A)+\fsl(B)$
modules are indexed in the obvious way. For the Segre case
write $T=U_{10}\op W_{10}$ and $N=U_{01}\ot W_{01}$. The remaining
relevant modules are as follows:
$$
\begin{aligned}
T&=V_{00001}\\
&=A_1\ot B_{100}\\
T^*&=V_{00010}\\
&=A_1\ot B_{001}\\
N&=N^*=V_{10000}\\
&=A_0\ot B_{010}\\
S^2T^*&=T^{*2}\op N\\
S^3T^*&=T^{*3}\op NT^*\\
S^3T^*\ot N&=(T^{*3}N\op T^{*2}T)\op ( N^2T^* \op \fg T^*\op NT\op T^*)\\
 T \ot N^*&= TN^*\op T^* \\
(N\ot N^*)^{\fr^c}\ot T^*&= N^2T^*\op NT\\
(T\ot T^*)^{\fr^c}\ot T^* &=T^{*2}T\op T_2T \op\fg T^*\op NT\op T^*
 \end{aligned}
$$
See \cite{LM1,LMseries} for an explanation of $T_2$.

The decompositions above show that there are $6$ components of $F_3$
on $\cf^1$ and four when we restrict to $\cf^n$. 

We may choose our model such that $\ee 1$ is
a highest weight vector (since it is in $\tbase\ii$). We may
also have $e^b$ be a highest weight vector for $N^*$.

Bertini easily kills the first component of $F_3$ as it has highest
weight vector $r^b_{111}$. In fact, as in the cases above, the
second Cartan component is also killed by Bertini. To see this note that  
 the two irreducible components of $S^3T^*$
respectively have highest weight vectors
$
r_{111}$ and a linear combination 
of $\ r_{11\overline 1}$ and $r_{12\overline 2}$.
The two Cartan components in $S^3T^*\ot N$ thus have highest weight
vectors
$
r_{111}^b$ and a linear combination of
  $r_{11\overline 1}^b$ and $r_{12\overline 2}^b
$. To see the second is zero, note that any linear combination of
$\ee 1,\ee 2$ is in $\tbase\ii$ and $\tsing q^b$. 

It remains to eliminate the second and fourth modules in $S^3T^*\ot N$, but neither
of these occurs in $(T\ot T^*)^{\fr^c}\ot T^* \op (N\ot N^*)^{\fr^c}\ot T^*$
and we are done with $F_3$.

The higher order invariants are safely left to the reader.
 \qed 

To compare with the $G(2,6)$ case in the standard model,
we have $\ee 1=\ee{(13)}$,$\ee {\overline 1}=\ee{(24)}$,
$q^a=q^{34}$,$q^b=q^{56}$ etc....

\end{document}